\documentclass{article}

\usepackage{amssymb, amsmath, amsthm}
\usepackage{standalone}
\usepackage{tikz}
\usepackage{multicol}
\usepackage{listings}
\usepackage{hyperref}
\hypersetup{hidelinks}

\newtheorem{theorem}{Theorem}[section]
\newtheorem{problem}[theorem]{Problem}
\newtheorem{lemma}[theorem]{Lemma}
\newtheorem{corollary}[theorem]{Corollary}
\theoremstyle{remark}
\newtheorem{remark}[theorem]{Remark}

\newcommand{\cell}[2]{\draw(#1,#2) --++(1,0) --++(0,1)--++(-1,0)--cycle;}
\newcommand{\one}[2]{
\filldraw[fill=white, draw=white] (#1+.1,#2+.1) rectangle (#1+.9,#2+.9);
}
\newcommand{\two}[2]{
\filldraw[fill=gray, draw=gray] (#1+.1,#2+.1) rectangle (#1+.9,#2+.9);
}
\newcommand{\three}[2]{
\filldraw[fill=black, draw=black] (#1+.1,#2+.1) rectangle (#1+.9,#2+.9);
}

\newcommand{\St}{\ensuremath{S_t}}
\newcommand{\Sr}{\ensuremath{S_r}}
\newcommand{\Srt}{\ensuremath{S_{r,t}}}

\newcommand{\sqmo}{\begin{tikzpicture}[scale=.25]
    \draw (0,0)--++(1,0)--++(0,1)--++(-1,0)--cycle;
    \draw (.5,0)--++(0,1);
    \draw (0,.5)--++(1,0);
    \end{tikzpicture}}

\newcommand{\tro}{\begin{tikzpicture}[scale=.25]
    \draw (.5,0)--++(0,1)--++(-.5,0)--++(0,-1)--++(1,0)--++(0,.5)--++(-1,0);
    \end{tikzpicture}}

\newcommand{\tmo}{\begin{tikzpicture}[scale=.25]
    \draw (1,0)--++(.5,0)--++(0,.5)--++(-1.5,0)--++(0,-.5)--++(1,0)--++(0,1)--++(-.5,0)--++(0,-1);
    \end{tikzpicture}}

\newcommand{\lmo}{\begin{tikzpicture}[scale=.25]
    \draw (0,.5)--++(1.5,0)--++(0,-.5)--++(-1.5,0)--++(0,1)--++(.5,0)--++(0,-1);
    \draw (1,0)--++(0,.5);
    \end{tikzpicture}}

\newcommand{\zmo}{\begin{tikzpicture}[scale=.25]
    \draw (1,0)--++(.5,0)--++(0,.5)--++(-1.5,0)--++(0,.5)--++(1,0)--++(0,-1)--++(-.5,0)--++(0,1);
    \end{tikzpicture}}

\newcommand{\clabel}[3]{\node (N) at (#1+.5,#2+.5){#3};}

\newcommand{\cells}[1]{
    \foreach \z in {0,...,#1}{
        \draw[-] (-\z,\z) --++(#1+1,0);
        \draw[-] (\z,-\z) --++(0,#1+1);
    }
}
\newcommand{\polycolored}[3]{
    \begin{tikzpicture}[scale=.39]
        \draw[-] #1;
        \clip #1;
        \cells{#2}
        #3 
    \end{tikzpicture}
}

\usepackage[backend=biber,
style=alphabetic]{biblatex}
\title{De Bruijn Polyominoes}
\author{D. M. Condon, Yuxin Wang, and E. Yang}
\date{Apr 2024}

\begin{document}

\maketitle

\begin{abstract}
    We introduce the notion of de Bruijn polyominoes, which generalizes the notions of de Bruijn sequences and arrays. Given a polyomino $p$ and a positive integer $n$, a $(p,n)$-de Bruijn polyomino is a colored polyomino $P$ with cells colored from $\{1,...,n\}$ such that every possible coloring of $p$ from $\{1,...,n\}$ occurs within $P$, not counting rotations. We discuss for some polyominoes $p$ and integers $n$ the shape of the smallest $(p,n)$-de Bruijn polyomino $P$, the construction of a coloring of $P$, and the enumeration of the de Bruijn colorings of $P$.
\end{abstract}

\section{Introduction and Background}

We are interested in polyominoes with colored cells.
A \textbf{polyomino} is a connected shape made from identical squares glued together edge-to-edge; these squares are called \textbf{cells}. We take these cells to be faces of the square lattice, with corners at integer Cartesian coordinates. The number of cells in a polyomino $P$ is called the \textbf{size} of $P$, denoted $|P|$. The polyominoes with up to four cells are depicted in Figure \ref{fig:intro_polyominoes}.

\begin{figure}
    \centering
    \begin{tikzpicture}[scale=.6]
    \cell 00
    \node (L1) at (.5,-.5) {monomino};
    
    \cell 30 \cell 40
    \node (L2) at (4,-.5) {domino};

    \cell {6.5}0 \cell {7.5}0 \cell {8.5}0
    \node (L3) at (8,-.5){straight tromino};

    \cell {14}0 \cell {11}0 \cell {12}0 \cell {13}0
    \node (L3) at (13,-.5){straight tetromino};

    \begin{scope}[yshift=-4cm,xshift=2.75cm]

    \cell {-4.5}{0} \cell {-3.5}{0}  \cell {-4.5}{1}
    \node (L4)[text width=2cm] at (-3,-.5) {L tromino};

    \cell {-1}0 \cell 01 \cell 00 \cell 10
    \node (L5)[text width=2cm] at (.5,-.5) {T tetromino};

    \cell 30 \cell 31 \cell 40 \cell 41
    \node (L7)[text width=2cm] at (4.5,-1) {square \ \ \ tetromino};

    \cell {6.5}1 \cell {7.5}{1} \cell {7.5}{0} \cell {8.5}{0}
    \node (L8) at (8,-.5) {Z tetromino};

    \cell {11}0 \cell {11}{1} \cell {12}{0} \cell {13}{0}
    \node (L9) at (12.5,-.5) {L tetromino};
    \end{scope}
\end{tikzpicture}
    \caption{The monomino, the domino, the trominoes, and the tetrominoes up to congruence.}
    \label{fig:intro_polyominoes}
\end{figure}

Polyominoes were popularized by Solomon Golomb in his book \textit{Polyominoes: Puzzles, Patterns, Problems, and Packings} \cite{Golomb}, which focused on tiling problems with polyominoes of certain sizes: the 1-cell \textbf{monomino}; the 2-cell \textbf{dominoes}; the 3-cell \textbf{trominoes}; the 4-cell \textbf{tetrominoes}, which appear in the popular game Tetris; and so on. 

Golomb considered congruent polyominoes to have the same shape, 
but in this paper we will only consider polyominoes to be the same {shape} up to translation. We call polyominoes that are the same shape \textbf{instances} of that shape. The orientations of the polyominoes in Figure \ref{fig:intro_polyominoes} are the orientations we will assume throughout this paper.

For any positive integer $n$, an \textbf{$n$-coloring} of a set of cells is a function $\chi$ that assigns each cell a color from $\{1,...,n\}$. A polyomino with $k$ cells has $n^k$ distinct $n$-colorings. Polyominoes of the same shape that have the same coloring up to translation are called \textbf{instances} of that coloring. 

Our original motivation is the following problem.

\begin{problem}\label{prob:sixteen_squares}
    Find the smallest 2-colored polyomino that includes exactly one instance of each of the 2-colorings of the square tetromino.
\end{problem}

The sixteen distinct 2-colorings of the square tetromino are depicted in Figure \ref{fig:sixteen_squares}, along with a 2-colored polyomino which includes one instance of each of these colorings but is not of minimum size.

\begin{figure}
    \begin{tikzpicture}[scale=.7]
\begin{scope}
    \foreach \z in {0,1,2,3,4,5,6}{
        \draw[-] (\z,0) --++(0,6);
        \draw[-] (0,\z) --++(6,0);
    }
    \foreach \x in {0,1,2,3,4,5}{
    \foreach \y in {0,1,2,3,4,5}{
        \two \x\y
    }
    }
    
    \draw[very thick, fill=white] (2,2) rectangle (4,4);
    \draw[very thick] (0,0)--++(6,0)--++(0,6)--++(-6,0)--cycle;

\one 01
\one 10
\one 13
\one 20
\one 21
\one 24
\one 25
\one 31
\one 34
\one 35
\one 41
\one 42
\one 44
\one 53
\one 54
\end{scope}

\begin{scope}[xshift=-10cm, scale=.7]
    \newcommand{\drawsquare}{
    \draw[very thick] (3*\x, 3*\y) --++(2,0)--++(0,2)--++(-2,0)--cycle;
    \draw[-] (3*\x+1, 3*\y) --++(0,2);
    \draw[-] (3*\x,3*\y+1) --++(2,0);
    }
    \foreach \x in {0,1,2,3}{
    \foreach \y in {0,1,2,3}{
        \drawsquare
    }
    }
    
    \foreach \x in {0,1,2,3}{
    \foreach \y in {0,1}{
        \two {3*\x}{3*\y+1}
    }
    }

    \foreach \x in {0,1,2,3}{
    \foreach \y in {0,2}{
        \two {3*\x+1}{3*\y+1}
    }
    }
    
    \foreach \x in {2,3}{
    \foreach \y in {0,1,2,3}{
        \two {3*\x}{3*\y}
    }
    }
    
    \foreach \x in {1,3}{
    \foreach \y in {0,1,2,3}{
        \two {3*\x+1}{3*\y}
    }
    }
\end{scope}
\end{tikzpicture}
        \caption{Left: The sixteen distinct 2-colorings of the square tetromino. White cells represent the color 1 and gray cells represent the color 2. Right: A polyomino with 32 cells, with a coloring that includes exactly one instance of each 2-coloring of the square tetromino. A smaller polyomino with this property exists, so this polyomino is $({\protect\sqmo},2)$-de Bruijn but does not solve Problem  \ref{prob:sixteen_squares}.}
    \label{fig:sixteen_squares}
\end{figure}

For a given polyomino $p$ and positive integer $n$ we define a \textbf{$(p,n)$-de Bruijn polyomino} to be an $n$-colored polyomino that includes exactly one instance of each $n$-coloring of $p$. We borrow the name \textit{de Bruijn} from de Bruijn sequences, which we discuss in Section \ref{sec:straight}. 
Problem \ref{prob:sixteen_squares} can be restated in these terms: find a $(\sqmo,2)$-de Bruijn polyomino of minimum size. Before proceeding, we invite the reader to try to solve Problem \ref{prob:sixteen_squares}.

Some natural questions that arise are: What shape should the polyomino be? Is this shape unique? Given an acceptable shape, how can we construct a de Bruijn coloring?  How many de Bruijn colorings are there?

For any positive integer $N$ and polyomino $p$, if $P$ is a polyomino of minimum size among those containing at least $N$ instances of $p$, we say that $P$ is \textbf{$(p,N)$-dense}. Our method for finding minimum de Bruijn polyominoes is to first construct dense polyominoes, and then search for satisfactory colorings of these. If a polyomino is $(p,n)$-de Bruijn and also $(p,n^{|p|})$-dense, we say it is \textbf{$(p,n)$-prismatic}.

In Section \ref{sec:connectivity} we show that the smallest set of cells that contains $N$ instances of any polyomino $p$ is necessarily connected, and is therefore a $(p,N)$-dense polyomino. In Section \ref{sec:straight} we show that if $p$ is a straight polyomino, then $(p,n)$-prismatic polyominoes are straight polyominoes which are colored by acyclic de Bruijn sequences. In Section \ref{sec:square}, we show that $(\sqmo,n)$-prismatic polyominoes are square polyominoes of side length $n^2+1$. In Section \ref{sec:tmo} we find the unique shape of $(\tmo,N^2)$-dense polyominoes. In Section \ref{sec:The_Tetrominoes}, we show how $(\sqmo,n)$-prismatic polyominoes can be transformed into $(\zmo,n)$-prismatic polyominoes, and how $(\tmo,N^2)$-dense polyominoes can be transformed into $(\lmo,N^2)$-dense polyominoes. In Section \ref{sec:Lmo}, we find the shapes of $(\tro,2)$-prismatic polyominoes.

\section{The Connectivity of Dense Polyominoes}\label{sec:connectivity}

In this section, we show that for any polyomino $p$ with at least two cells, and any positive integer $N$, if $S$ is a set of cells of minimum size containing at least $N$ instances of $p$, then $S$ must be connected, and the instances themselves must be connected through shared cells.

\begin{lemma}\label{thm:connected}
    For any positive integer $N$ and any polyomino $p$ with at least two cells, if $S$ is a set of cells of minimum size that contains at least $N$ instances of $p$ then $S$ is connected. In particular, $S$ is a $(p,N)$-dense polyomino.
\end{lemma}

\begin{proof}

    We proceed by contradiction. Let $S$ be a set of cells of minimum size that contains at least $N$ instances of $p$, and suppose that $S$ is disconnected.
    Since $p$ has at least two cells, it occupies at least two rows or at least two columns. We will assume without a loss of generality that $p$ occupies at least two columns.

    Consider the case that $S$ has only two connected components, $C_1$ and $C_2$. 
    Suppose each $C_i$ contains $N_i$ instances of $p$, so that $N_1+N_2 \geq N$. We construct the set $S'$ by translating $C_2$ 
    to a set $C_2'$ of the same shape
    such that some left-most cell of $C_2'$ is also a right-most cell of $C_1$. Note that $S' = C_1 \cup C_2'$ is connected, and has fewer cells than $S$. Our task is now to show that $S'$ contains at least $N$ instances of $p$. 
    Since $C_1$ and $C_2'$ only overlap in one column, and $p$ spans multiple columns, no instance of $p$ is shared by $C_1$ and $C_2'$.
    Since $C_1$ contains $N_1$ instances of $p$, and $C_2'$ contains $N_2$ instances of $p$,  $S'$ contains at least $N_1+N_2 \geq N$ instances of $p$. Since $S'$ has fewer cells than $S$, we have arrived at a contradiction.

    In the case that $S$ has more than two components, we may apply the same argument with two of the components to arrive at the same contradiction.
\end{proof}

Dense polyominoes satisfy an even stronger notion of connectedness, as described by the following lemma.

\begin{lemma}\label{thm:strongly_connected}
    For any positive integer $N$, any polyomino $p$ with at least two cells, and any $(p,N)$-dense polyomino $P$, let $G$ be the graph whose vertices are the instances of $p$ in $P$ with adjacency between those instances that share at least one cell. Then $G$ is a connected graph.
\end{lemma}

\begin{proof}
    We proceed by contradiction. Suppose $G$ is not connected, and has multiple components $C_1,...,C_m$. Then for $i \neq j$, the cells belonging to instances in $C_i$ cannot be shared with cells belonging to instances in $C_j$. We can construct a new set of cells $P'$ from $P$ by translating all of the cells belonging to instances in $C_1$ far away from the rest of the cells so that $P'$ is disconnected; this construction does not remove any instances of $p$ from the shape, nor change the number of cells. Then $P'$ is a set of cells of minimum size that contains at least $N$ instances of $p$, but is disconnected, so we have arrived at a contradiction by Lemma \ref{thm:connected}.
\end{proof}

\section{Straight Polyominoes}\label{sec:straight}

A \textbf{straight polyomino} is a polyomino with at least two cells, with all cells occupying a single row of the square lattice.\footnote{A polyomino with all cells occupying a single column is also usually considered a straight polyomino, but as with other polyominoes we are restricting our focus to a single orientation.} Straight polyominoes of the same size are necessarily the same shape.

For $k$ a positive integer and $A$ a finite alphabet, a \textbf{de Bruijn sequence} of order $k$ over $A$ is a cyclic sequence which includes each sequence of length $k$ over $A$ exactly once as a subsequence. For example, the sequence $$(0,0,1,1,1,0,1,0)$$
includes each sequence of length 3 over $\{0,1\}$ when interpreted cyclically.

In this section, we will show when $p$ is a straight polyomino that $(p,n)$-prismatic polyominoes exist, are straight polyominoes, and have de Bruijn colorings that come from de Bruijn sequences.

\begin{lemma}\label{thm:straightpolyominoes}
    For any positive integers $N$ and $k$, with $k>1$, let $p$ be the straight polyomino of size $k$. Then the unique shape of a $(p,N)$-dense polyomino is the straight polyomino of size $N+k-1$.
\end{lemma}

\begin{proof}   
    Let $P$ be a $(p,N)$-dense polyomino.
    Since straight polyominoes only share cells if they occupy the same row of the square lattice, it follows from Lemma \ref{thm:strongly_connected} that all instances of $p$ in $P$ must occupy one row. Since $P$ is $(p,N)$-dense it has no cells which do not belong to an instance of $p$, so every cell of $P$ must occupy the same row; therefore $P$ is a straight polyomino. 
    The number of instances of $p$ within $P$ is then $|P| - k+1$, so $P$ must be a straight polyomino of size $N+k-1$.
\end{proof}

We will now show that for $p$ a straight polyomino of size $k > 1$, the straight polyomino $P$ of size $n^k+k-1$ admits a $(p,n)$-de Bruijn coloring.

We will represent $n$-colorings of $P$ as finite sequences over $\{1,...,n\}$; the $i$th term in the sequence represents the color of the $i$th cell from the left. The colorings of instances of $p$ correspond to subsequences of length $k$. The $(p,n)$-de Bruijn colorings of $P$ therefore correspond to those sequences of length $n^k+k-1$ over $\{1,...,n\}$ that include each possible subsequence of length $k$ exactly once. These are essentially acyclic de Bruijn sequences.

We can construct acyclic de Bruijn sequences of order $k$ over $\{1,...,n\}$ from cyclic de Bruijn sequences: choose any starting point in a cyclic de Bruijn sequence 
and read the sequence of length $n^k+k-1$ from that point, so that the first $k-1$ terms are identical to the last $k-1$ terms. 
See Figure \ref{fig:bruijn_sequences} for an example. In fact, any acyclic de Bruijn sequence can be constructed in this way.

\begin{figure}
    \centering
    \newcommand{\onelabel}[2]{
\filldraw[fill=white, draw=white] (#1+.1,#2+.1) rectangle (#1+.9,#2+.9);
\node (O#1#2) at (#1+.5,#2+.5) {1};
}
\newcommand{\twolabel}[2]{
\filldraw[fill=white, draw=white] (#1+.1,#2+.1) rectangle (#1+.9,#2+.9);
\node (O#1#2) at (#1+.5,#2+.5) {2};
}

\begin{tikzpicture}[scale=.6]
        \draw [thick] (1,0)--++(19,0)--++(0,1)--++(-19,0)--cycle;
        \foreach \x in {1,...,19}{
            \draw[thick] (\x,0)--(\x, 1);
            \one \x 0 
        }

        \foreach \x in {9,11,13,14,15,16,21,24,25,27}{
            \two {\x-8} 0
        }
        
        \begin{scope}[yshift=2cm]

        \foreach \x in {2,...,19}{
            \node (C\x) at (\x,.25) {,};
            \onelabel \x 0 
        }

        \foreach \x in {9,11,13,14,15,16,21,24,25,27}{
            \twolabel {\x-8} 0
        }
        \end{scope}
        
        \begin{scope}[yshift=4cm]

        \foreach \x in {2,...,16}{
            \node (C\x) at (\x,.25) {,};
            \onelabel \x 0 
        }

        \foreach \x in {9,11,13,14,15,16,21,24}{
            \twolabel {\x-8} 0
        }
        \node (Op) at (1,.5) {(};
        \node (Cp) at (17,.5) {)};
        \end{scope}
\end{tikzpicture}
    \caption{
    Top: a cyclic de Bruijn sequence which which includes every sequence over $\{1,2\}$ of length 4 exactly once. 
    Middle: an acyclic de Bruijn sequence constructed from the cyclic de Bruijn sequence by reading the sequence of length 19 starting from the first entry.
    Bottom: a $(p,2)$-prismatic polyomino colored by the acyclic de Bruijn sequence, where $p$ is the straight polyomino of size 4. }
    \label{fig:bruijn_sequences}
\end{figure}

T. van Aardenne-Ehrenfest and N. G. de Brujin \cite{BE} showed that the number of cyclic de Bruijn sequences of order $k$ over an alphabet of size $n$ is 
$$\frac{(n!)^{n^{k-1}}}{n^k}.$$
These cyclic de Bruijn sequences each have $n^k$ entries, since each entry is the start of one substring of length $k$, so there are $n^k$ acyclic de Bruijn sequences corresponding to each cyclic de Bruijn sequence. The number of acyclic de Bruijn sequences of order $k$ over $\{1,...,n\}$ is then
$$(n!)^{n^{k-1}}.$$
This formula, together with Lemma \ref{thm:straightpolyominoes}, gives us the following theorem.

\begin{theorem}\label{thm:sequence_count}
    For any positive integers $n$ and $k$, with $k>1$, let $p$ be the straight polyomino of size $k$. Then there are $(n!)^{n^{k-1}}$ distinct $(p,n)$-prismatic polyominoes, which are all straight polyominoes of size $n^k+k-1$.
\end{theorem}

\section{The Square Tetromino} \label{sec:square}

In this section we will show that $(\sqmo,n)$-prismatic polyominoes exist, are square polyominoes, and have colorings that correspond to certain de Bruijn arrays:
For positive integers $k_1$ and $k_2$, and a finite alphabet $A$ of size $n$, a \textbf{de Bruijn array} is a matrix with entries in $A$ that contains as a connected submatrix exactly one copy of each possible $k_1 \times k_2$ matrix with entries in $A$.

De Bruijn arrays are less well studied than de Bruijn sequences; for example, little is known about the number of de Bruijn arrays for general parameters. Some small cases have been calculated, and Mark Dow compiled some of these results \cite{Dow}.
A means of constructing certain examples of these arrays was given by J. C. Cock \cite{Cock}, which we demonstrate in Section \ref{sec:cock}.  

\subsection{The Shape of Dense Polyominoes}

We will use a generalization of Pick's Theorem due to Dale E. Varberg \cite{DV}.

\begin{theorem}[Generalized Pick's Theorem]
    For any integer polygon with area $A$, $B$ boundary points, $I$ interior points, and Euler characteristic $E$,
    $$\frac B2 = A -I+E.$$
\end{theorem}

Here, an \textbf{integer polygon} is a simple polygon on the Cartesian plane with vertices at integer coordinates; the shape may have holes, the boundaries of which must be simple polygons with vertices at integer coordinates; edges of the boundary may only meet at integer coordinates.
Points with integer coordinates along the boundary are called \textbf{boundary points}, and points with integer coordinates within the interior are called \textbf{interior points}. 

Observe that a polyomino is an integer polygon. Its area is its size. Its interior points correspond to the instances of the square tetromino contained by the polyomino. Its Euler characteristic is straightforward to compute: $E$ is 1 minus the number of holes. In particular, $E \leq 1$. Applying Pick's Theorem to polyominoes, we arrive at the following corollary:

\begin{corollary}\label{thm:pick}
    For any polyomino with size $A$, $B$ boundary points, and containing $I$ instances of the square tetromino,
    $$\frac B2 \leq A -I+1.$$
\end{corollary}

We are now ready to prove that $(\sqmo,N^2)$-dense polyominoes are square polyominoes.

\begin{lemma}\label{thm:big_square}
    For any positive integer $N$, 
    the unique shape of a $(\sqmo,N^2)$-dense polyomino is the square polyomino of side length $N+1$.
\end{lemma}

Our proof makes use of the \textbf{width} (resp. \textbf{height}) of a polyomino, which is the number of columns (resp. rows) its cells occupy.

\begin{proof}     
    Note that the square polyomino of side length $N+1$ has $(N+1)^2$ cells and contains $N^2$ instances of $\sqmo$.
    
    Let $P$ be any polyomino having at most $(N+1)^2$ cells and containing at least $N^2$ instances of $\sqmo$; then $A \leq (N+1)^2$ and $ I \geq N^2$
    in the language of Corollary \ref{thm:pick} applied to $P$.
    We will show that $P$ must be the square polyomino of side length $N+1$.
    
    Suppose that $P$ has width $w$ and height $h$. Then the outer boundary of $P$ is a simple polygon with length at least $2w+2h$, so $2w+2h \leq B$ in the language of Corollary~\ref{thm:pick} applied to $P$.  
    It follows from Corollary~\ref{thm:pick} that 
    $$w + h \leq \frac B2 \leq A-I+1 \leq (N+1)^2 - N^2 +1= 2N+2.$$

    Our polyomino $P$ fits inside a rectangular polyomino $R$ with width $w$ and height $h$.
    The number of instances of $\sqmo$ in $R$ is $(w-1)(h-1)$. Since $w + h \leq 2N+2$, 
    it is straightforward to check with calculus that $(w-1)(h-1) \leq N^2$ with equality when $w = h = N+1$. 
    Since $R$ must contain at least $N^2$ instances of $\sqmo$, $R$ is then a square polyomino of side length $N+1$ and contains exactly $N^2$ instances of~$\sqmo$.
    Since every cell of $R$ belongs to at least one instance of $\sqmo$, any polyomino that is smaller than $R$ and fits inside $R$ contains fewer than $N^2$ instances of $\sqmo$. Since $P$ contains at least $N^2$ instances of $\sqmo$, it must be that $P = R$, so $P$ is a square polyomino of side length $N+1$.
\end{proof}

\subsection{A Coloring Construction} \label{sec:cock}

Lemma \ref{thm:big_square} implies that if $(\sqmo,n)$-prismatic polyominoes exist then they are squares of side length \mbox{$n^2+1$}. Now we will show that these polyominoes admit $(\sqmo,n)$-de Bruijn colorings.

\begin{theorem}\label{thm:prismatic_square}
    For any positive integer $n$, there are at least $n!^{n} \cdot (n^2)$ distinct $(\sqmo,n)$-prismatic polyominoes, which are all square polyominoes of side length $n^2+1$.
\end{theorem}

Our proof adapts an algorithm for constructing de Bruijn arrays which is due to Cock \cite{Cock}. An example of the construction is given in Figure \ref{fig:cock}.

\begin{proof}
Let a positive integer $n$ be given.
By Lemma \ref{thm:big_square}, it suffices to show that $P$ the square polyomino of side length $n^2+1$ admits a $(\sqmo,n)$-de Bruijn coloring.
We will construct an $n$-coloring for $P$, and then demonstrate that each $n$-coloring of $\sqmo$ occurs at least once (and therefore exactly once) in $P$. Lastly, we count the number of different de Bruijn colorings given by the construction.

The construction uses the following algorithm:
\begin{enumerate}
    \item Choose an arbitrary de Bruijn sequence $r_1$ of order $2$ over $\{1,...,n\}$, with $r_1 =  (x_1,...,x_{n^2}).$
    \item Choose an arbitrary permutation $\sigma$ of $\{1,...,n^2\}$, with $\sigma = (\sigma_1,...,\sigma_{n^2})$.
    \item We define $r_2,...,r_{n^2+1}$ to be cyclic permutations of $r_1$, where $r_{i+1}$ is the cyclic permutation of $r_{i}$  by $\sigma_i$, for $i \in \{1,2,...,n^2\}$.
    \item We index the rows and columns of $P$ from top to bottom and left to right, starting from 1. We color the $j$th cell in the $i$th row of $P$ by the $j$th entry of $r_i$, for $i = 1,..., n^2+1$ and $j = 1,..., n^2$. We make the final cell in each row the same color as the first cell in that row.
\end{enumerate}

Given a coloring of $\sqmo$ with the top row colored $(w,x)$ and bottom row colored $(y,z)$, we can identify a place where this coloring of $\sqmo$ occurs within the coloring of $P$ as follows:
\begin{enumerate}
    \item First, determine the positions of $(w,x)$ and $(y,z)$ as subsequences of $r_1$, and find $k$ such that $(y,z)$ occurs $k$ positions after $(w,x)$ within $r_1$ (interpreted cyclically).
    \item Find the unique $i$ such that $\sigma_{i}=k$.
    \item Find the index $j$ of the unique $w$ within $r_i$ that is followed immediately by an $x$.
    \item Note that the $j$th and $j+1$st entries of $r_{i+1}$ are $y$ and $z$ respectively, since $r_{i+1}$ is the cyclic permutation of $r_i$ by $\sigma_i = k$.
    \item Then the instance of $\sqmo$ that has its top-left cell in row $i$ and column $j$ of $P$ has its top row colored $(w,x)$ and its bottom row colored $(y,z)$. 
\end{enumerate}
Thus, each coloring of $\sqmo$ occurs within the coloring of $P$, as was to be shown.

Finally, we note that the construction depends on a choice of
$r_1$ (a de Bruijn sequence of order $2$ over an alphabet of size $n$), a particular starting point in $r_1$, and $\sigma$ an arbitrary permutation of $\{1,2,...,n^2\}$, and that different choices for these will result in different colorings. The construction therefore gives
$n!^{n} \cdot (n^2)!$ distinct $(\sqmo,n)$-de Bruijn colorings.
\end{proof}

\begin{figure}
    \centering
    \setlength{\tabcolsep}{2pt}
\renewcommand{\arraystretch}{1}
\begin{center}
    \begin{tabular}[b]{rccccccccc}
        $r_1=$ & (1,& 1,& 2,& 2,& 3, &3, & 1, & 3, &2)\\
        $r_2=$ & (1, & 2,& 2,& 3, &3, & 1, & 3, &2, &1)\\
        $r_3=$ & (2,& 3, &3, & 1, & 3, &2, &1, &1, &2)\\
        $r_4=$ & (1, & 3, &2, &1, &1, &2, &2,& 3, &3)\\
        $r_5=$ & (1, &2, &2,& 3, &3, &1, & 3, &2, &1)\\
        $r_6=$ & (1, & 3, &2, &1, & 1, &2, &2,& 3, &3)\\
        $r_7=$ & (2,& 3, &3, &1, & 3, &2, &1, & 1, &2)\\
        $r_8=$ & (1, &2, &2,& 3, &3, &1, & 3, &2, &1)\\
        $r_9=$ & (1, & 1, &2, &2,& 3, &3, &1, & 3, &2)\\
        $r_{10}=$ & (1, & 1, &2, &2,& 3, &3, &1, & 3, &2)\\
    \end{tabular}
    \hspace{1cm}
    \begin{tikzpicture}[scale=.4, y=-1cm]
    \draw[very thick] (1,1)--++(10,0)--++(0,10)--++(-10,0)--cycle;
    \foreach \z in {2,...,10}{
        \draw[-] (\z,1) --++(0,10);
        \draw[-] (1,\z) --++(10,0);
    }

\one{1}{10}
\one{2}{10}
\two{3}{10}
\two{4}{10}
\three{5}{10}
\three{6}{10}
\one{7}{10}
\three{8}{10}
\two{9}{10}
\one{10}{10}
\one{1}{9}
\one{2}{9}
\two{3}{9}
\two{4}{9}
\three{5}{9}
\three{6}{9}
\one{7}{9}
\three{8}{9}
\two{9}{9}
\one{10}{9}
\one{1}{8}
\two{2}{8}
\two{3}{8}
\three{4}{8}
\three{5}{8}
\one{6}{8}
\three{7}{8}
\two{8}{8}
\one{9}{8}
\one{10}{8}
\two{1}{7}
\three{2}{7}
\three{3}{7}
\one{4}{7}
\three{5}{7}
\two{6}{7}
\one{7}{7}
\one{8}{7}
\two{9}{7}
\two{10}{7}
\one{1}{6}
\three{2}{6}
\two{3}{6}
\one{4}{6}
\one{5}{6}
\two{6}{6}
\two{7}{6}
\three{8}{6}
\three{9}{6}
\one{10}{6}
\one{1}{5}
\two{2}{5}
\two{3}{5}
\three{4}{5}
\three{5}{5}
\one{6}{5}
\three{7}{5}
\two{8}{5}
\one{9}{5}
\one{10}{5}
\one{1}{4}
\three{2}{4}
\two{3}{4}
\one{4}{4}
\one{5}{4}
\two{6}{4}
\two{7}{4}
\three{8}{4}
\three{9}{4}
\one{10}{4}
\two{1}{3}
\three{2}{3}
\three{3}{3}
\one{4}{3}
\three{5}{3}
\two{6}{3}
\one{7}{3}
\one{8}{3}
\two{9}{3}
\two{10}{3}
\one{1}{2}
\two{2}{2}
\two{3}{2}
\three{4}{2}
\three{5}{2}
\one{6}{2}
\three{7}{2}
\two{8}{2}
\one{9}{2}
\one{10}{2}
\one{1}{1}
\one{2}{1}
\two{3}{1}
\two{4}{1}
\three{5}{1}
\three{6}{1}
\one{7}{1}
\three{8}{1}
\two{9}{1}
\one{10}{1}
    
\end{tikzpicture}
\end{center}
\setlength{\tabcolsep}{6pt}
    \caption{We demonstrate the construction of a $({\protect\sqmo},3)$-de Bruijn coloring of the square polyomino of side length 10. Let $r_1 = (1,1,2,2,3,3,1,3,2)$. Let $\sigma = (1,2,3,4,5,6,7,8,9)$. Then $r_{i+1}$ is a cyclic permutation of $r_{i}$ by $i$. These cyclic sequences are listed above on the left. The corresponding prismatic polyomino is on the right, with 1 represented by white, 2 represented by gray, and 3 represented by black. The last cell in each row is colored identically to the first cell in that row.
    To find the instance of the square tetromino with top row colored $(2,1)$ and bottom row colored $(2,3)$, we first note that $(2,3)$ appears 4 positions after $(2,1)$ in $r_1$, interpreted cyclically. Since $\sigma_{i}=4$ when $i=4$, we note the position of the subsequence $(2,1)$ in $r_4$: $2$ and $1$ are the third and fourth entries of $r_4$. Thus, $2$ and $3$ are the third and fourth entries of $r_5$. The instance in question then spans the fourth and fifth rows, and the third and fourth columns.
    }
    \label{fig:cock}
\end{figure}

\begin{remark}
    Cock's construction generates 96 $(\sqmo,2)$-prismatic polyominoes, and 78,382,080 $(\sqmo,3)$-prismatic polyominoes.     
    There are 800 $(\sqmo,2)$-prismatic polyominoes in total, so many of these are not generated by Cock's construction. Appendix \ref{app:16sols} includes code which generates all 800 of these.
\end{remark}

\section{The T Tetromino}\label{sec:tmo}

In this section we discuss $(\tmo,N^2)$-dense polyominoes. We will refer to the cells in the T tetromino as the \textit{left}, \textit{central}, \textit{right}, and \textit{top} cell as shown in Figure~\ref{fig:ziggy_stardust}.

\begin{figure}
\begin{center}
\begin{tikzpicture}[scale=1.25]
    \draw[-] (1,0)--++(0,2)--++(1,0)--++(0,-2)--++(-2,0)--++(0,1)--++(3,0)--++(0,-1)--++(-1,0);
    \node (L) at (.5,.5){left};
    \node (C) at (1.5,.5){central};
    \node (R) at (2.5,.5){right};
    \node (T) at (1.5,1.5){top};
\end{tikzpicture}
\hfill
\begin{tikzpicture}[scale=.5]
    \foreach \x in {0,1,2,3,4}{
        \draw (\x,0) rectangle (9 - \x, 1+\x);
    }
\end{tikzpicture}
\hfill
\end{center}
\caption{Left: The cells of the T tetromino. Right: A ziggurat of height 5.}
\label{fig:ziggy_stardust}
\end{figure}

We define a \textbf{ziggurat of height $N$} to be a polyomino with $N$ center-aligned rows having odd lengths $1,3,...,2N-1$ from top to bottom. See Figure \ref{fig:ziggy_stardust}. The shape has $N^2$ cells in total. The cells above the bottom row are the top cells of the instances of the T tetromino within the shape, and so the ziggurat of height $N$ contains $(N-1)^2$ instances of the T tetromino. Our main result in this section is the following theorem.   

\begin{theorem}\label{thm:ziggy_stardust}
    For any positive integer $N$, the unique shape of a $(\tmo,N^2)$-dense polyomino is the ziggurat of height $N+1$. 
\end{theorem}

\subsection{The Dimensions of Dense Polyominoes}
We start with a technical lemma.

\begin{lemma}\label{thm:thick_rows}
    For any positive integer $N$, let $P$ be a $(\tmo,N)$-dense polyomino. Then no row of $P$ below the top row can have fewer than three cells.
\end{lemma}

\begin{proof}
    For $N$ a positive integer, let $P$ be a $(\tmo,N)$-dense polyomino. Note that each cell in $P$ must belong to an instance of $\tmo$. Suppose some row $r$ of $P$ has fewer than three cells; we will show that $r$ is the top row of $P$. 

    Since $r$ has fewer than three cells, no instance of $\tmo$ in $P$ has cells in both $r$ and the row directly above $r$. In fact, no instance of $\tmo$ with a cell in $r$ or any row below $r$ can have a cell in any row above $r$; so every instance of $\tmo$ in $P$ either has all cells in rows above $r$, or has all cells in $r$ and rows below $r$.
    Any instances of $\tmo$ with all cells above $r$ are then disconnected, in the sense of Lemma \ref{thm:strongly_connected}, from those instances of $\tmo$ with cells in $r$ and rows below $r$; therefore, there are no instances of $\tmo$ with cells above $r$. Then $P$ has no cells above $r$, so $r$ is the top row of~$P$.
\end{proof}

    The following lemma bounds the width and height of $(\tmo,N^2)$-dense polyominoes.

    \begin{lemma}\label{thm:ziggurat_box}
        For any positive integer $N$, if $P$ is a $(\tmo,N^2)$-dense polyomino, then $P$ has height at most $N+1$ and width at most $2N+1$.
    \end{lemma}

    \begin{proof}
    Suppose that for some positive integer $N$, $P$ is a $(\tmo,N^2)$-dense polyomino with width~$w$ and height~$h$.

    At least $N^2$ cells in $P$  are central cells for some instance of $\tmo$.
    The top row of $P$ contains at least one cell, which cannot be a central cell because it is in the top row.
    By Lemma \ref{thm:thick_rows}, each row below the top row has a distinct left-most and right-most cell.
    There are $h-1$ such left-most cells and $h-1$ such right-most cells, and none of these are central cells for any instance of $\tmo$.
    Therefore, $P$ has at least $N^2+2(h-1)+1 = N^2+2h-1$ cells.
    Since the ziggurat of height $N+1$ with $(N+1)^2$ cells gives an upper bound on the size of $P$, 
    $$N^2 +2h-1 \leq |P| \leq N^2 + 2N+1$$
    and so
    $h \leq N+1.$

    Similarly, the top cell in any column is not the central cell for any instance of $\tmo$. There are $w$ columns. Therefore
    $$N^2+w \leq |P| \leq N^2+2N+1$$
    and so 
    $w \leq 2N+1.$
    \end{proof}

    This lemma implies that any $(\tmo,N^2)$-dense polyomino fits inside a rectangular polyomino of height $N+1$ and width $2N+1$.

    \subsection{A Proof by Linear Program}
    We are now ready to prove Theorem \ref{thm:ziggy_stardust}.
    \begin{proof}
    Suppose that for some positive integer $N$, $P$ is a $(\tmo,N^2)$-dense polyomino. By Lemma \ref{thm:ziggurat_box}, $P$ fits inside a rectangular polyomino $R$ of height $N+1$ and width $2N+1$. We may assume without a loss of generality that some cell of $P$ occupies the top row of $R$.
    Indexing the rows of $R$ from top to bottom, let $r_i$ be the number of cells in the $i$th row of $R$ that belong to $P$, for $i=1,...,N+1$. Let $a_i$ be the number of cells in the $i$th row of $R$ that are the top cell of some instance of $\tmo$ within $P$, for $i=1,...,N$. Then
    $$0 \leq a_i \leq r_i \leq 2N+1.$$
    As the ziggurat of height $N$ provides an upper bound on the size of $P$,
    \begin{equation}\label{eq:inq2}
        \sum_{i=1}^{N+1} r_i \leq (N+1)^2, \quad \text{ for } i = 1,...,N.
    \end{equation}
    Since $P$ contains at least $N^2$ instances of $\tmo$, each having one top cell,
    \begin{equation}\label{eq:inq3}
        N^2 \leq \sum_{i=1}^{N} a_i.
    \end{equation}
    For each top cell in row $i$, there is a central cell directly below it in row $i+1$, which cannot be the left-most or right-most cell in row $i+1$. Therefore,
    \begin{equation}\label{eq:inq1}
        a_i+2 \leq r_{i+1}, \quad \mbox{ for } i=1,...,N.
    \end{equation}
    Combining these weak inequalities, 
\begin{align*}
    (N+1)^2 &= 1+2N+N^2\\
    \mbox{by } \eqref{eq:inq3} \quad &\leq 1 + 2N + \sum_{i=1}^N a_i\\
    & = 1 + \sum_{i=1}^N (a_i+2)\\
    &\leq r_1 + \sum_{i=1}^N (a_i+2)\\
    \mbox{by } \eqref{eq:inq1} \quad &\leq \sum_{i=1}^{N+1} r_i\\
    \mbox{by } \eqref{eq:inq2} \quad &\leq (N+1)^2.
\end{align*}
Since each of these weak inequalities must occur with equality, we find
$$\sum_{i=1}^N a_i = N^2, \quad r_1 = 1, \quad  a_i+2 = r_{i+1}, \quad \mbox{ and } \quad \sum_{i=1}^{N+1} r_i = (N+1)^2.$$
Since
\begin{equation*}
    r_{i+1} = a_i+2 \leq r_i+2,  \quad \mbox{for } i=1,...,N,
\end{equation*}
it follows that each $r_i \leq 2i-1$. We then have the inequality
$$(N+1)^2 = \sum_{i=1}^{N+1} r_i \leq \sum_{i=1}^{N+1}2i-1 = (N+1)^2.$$
This inequality must occur with equality, so
\begin{equation*}
    r_{i} = 2i-1   \quad \mbox{for } i=1,...,N+1, \quad \mbox{ and } \quad  a_i = 2i-1   \quad \mbox{for } i=1,...,N.
\end{equation*}

A straightforward argument by induction on $j$ shows that the first $j$ rows of $P$ form a ziggurat of height $j$, so that $P$ itself is a ziggurat of height $N+1$.
    \end{proof}

    \begin{remark}
        Theorem \ref{thm:ziggy_stardust} shows that if a $(\tmo,n)$-prismatic polyomino exists, it must be a ziggurat of height $n^2+1$. It remains to show that this ziggurat admits a $(\tmo,n)$-de Bruijn coloring.

        We exhibit such a coloring for $n=2$ in Figure \ref{fig:T_prismatic}. A branch-and-prune search based on the code in Appendix \ref{app:16sols} shows there are 168 distinct $(\tmo, 2)$-prismatic polyominoes.
    \end{remark}

\section{The Z and L Tetrominoes}\label{sec:The_Tetrominoes}
This section discusses a polyomino transformation which gives a bijection between $(\sqmo,n)$-prismatic and $(\zmo,n)$-prismatic polyominoes, as well as between $(\tmo,N^2)$-dense and $(\lmo,N^2)$-dense polyominoes.

We associate each cell on the square lattice with the Cartesian coordinates of its bottom-left corner.
Let $f: \mathbb Z^2 \rightarrow \mathbb Z^2$ be the function $f:(x,y) \mapsto (x-y,y)$, a bijection on the set of cells in the square lattice which shifts each row one cell to the right relative to the row above it.

For $S$ a set of cells, we define
$f(S): = \{f(c): c \in S\}.$
If $S$ has an associated coloring $\chi$, then $f$ and $\chi$ induce a coloring of $f(S)$, namely the coloring $\chi'$ where $\chi'(f(c)) = \chi(c)$ for each $c \in S$. So $f$ induces a transformation of (possibly colored) sets of cells. We call this transformation the \textbf{row shift}, and denote it with the letter $f$. We note that the row shift is invertible.

The row shift transforms a polyomino by shifting each row one cell to the right relative to the row above, resulting in another set of cells which is a polyomino so long as it is connected. For example, the row shift of a square tetromino is always a Z tetromino. The row shift of a Z tetromino is a disconnected set of cells.

\subsection{A Bijection for Square and Z Tetrominoes}\label{sec:sq_and_z}
In this section we will show how to construct a $(\zmo, n)$-prismatic polyomino from a $(\sqmo, n)$-prismatic polyomino. We start with some technical lemmas.

\begin{lemma}\label{thm:sq_z_bij}
    Let $p$ and $q$ be sets of cells such that $f(p) = q$. Then $p$ is an instance of the square tetromino if and only if $q$ is an instance of the Z tetromino.
\end{lemma}

\begin{proof}
    Let $p$ and $q$ be sets of cells such that $f(p)=q$. Note that $p$ is an instance of the square tetromino if and only if it has four cells that have coordinates $(m,n)$, $(m+1,n)$, $(m,n-1)$ and $(m+1,n-1)$ for some integers $m$ and $n$. Similarly, $q$ is an instance of the Z tetromino if and only if it has four cells that have coordinates $(x,y)$, $(x+1,y)$, $(x+1,y-1)$ and $(x+2,y-1)$ for some integers $x$ and $y$.
    
    Suppose $p$ is a square tetromino, whose cells have coordinates $(m,n)$, $(m+1,n)$, $(m,n-1)$ and $(m+1,n-1)$ for some integers $m$ and $n$. Then $q$ is
    \begin{align*}
        f(p) = \{(m-n,n), (m+1-n,n), (m-n-1,n-1), (m+1-n-1,n-1)\}.
    \end{align*}
    Let $x = m - n$ and $y = n$. Then the coordinates for $q$ are 
    \begin{align*}
        \{(x,y), (x + 1, y), (x + 1,y-1), (x + 2, y-1)\}
    \end{align*}
    so $q$ is a Z tetromino. The other direction can be shown in the same way.
\end{proof}

\begin{lemma}\label{thm:dense_bij}
    Let $P$ and $Q$ be sets of cells such that $f(P)=Q$. Then for any positive integer $N$, $P$ is a $(\sqmo,N)$-dense polyomino if and only if $Q$ is a $(\zmo,N)$-dense polyomino.
\end{lemma}

\begin{proof}
    It follows from Lemma \ref{thm:sq_z_bij} that the row shift $f$ induces a bijection between the instances of $\sqmo$ in $P$ and the instances of $\zmo$ in $Q$. 
    If $P$ is a $(\sqmo,N)$-dense polyomino, then $Q$ contains at least $N$ instances of $\zmo$ and has size $|P|$. If any smaller polyomino $Q'$ had at least $N$ instances of $\zmo$, then $f^{-1}(Q)$ would be a set of cells smaller than $P$ which contained $N$ instances of $\sqmo$; then $P$ would not be $(\sqmo,N)$-dense, a contradiction. So no such $Q'$ exists, thus $Q$ must be $(\zmo,N)$-dense.  The reverse direction follows by identical reasoning.
\end{proof}

\begin{lemma}\label{thm:color_bij}
    For any positive integer $n$, a set of cells $S$ with an $n$-coloring includes exactly one instance of each $n$-coloring of the square tetromino if and only if $f(S)$ includes exactly one instance of each $n$-coloring of the Z tetromino. 
\end{lemma}

\begin{proof}
    It follows from Lemma \ref{thm:sq_z_bij} that the row shift $f$ induces a bijection between the instances of $\sqmo$ in any set of cells $S$ and the instances of $\zmo$ in $f(S)$. 
    
    It is straightforward to check that if $p$ and $p'$ are colored sets of cells with the same shape, then $p$ and $p'$ have the same coloring if and only if $f(p)$ and $f(p')$ have the same coloring.

    Suppose that $S$ is a set of cells with an $n$-coloring that includes exactly one instance of each $n$-coloring of $\sqmo$.
    Then $S$ includes $n^4$ distinctly colored instances of $\sqmo$. The row shift must then map each of these to distinctly colored instances of $\zmo$ in $f(S)$, so that $f(S)$ includes exactly one instance of each $n$-coloring of $\zmo$. The reverse direction uses identical reasoning.
\end{proof}

We are now ready to prove our main theorem for this section.

\begin{figure}
    \centering
    \begin{tikzpicture}[scale=.5]
    \begin{scope}[xshift=-10cm]
    \draw [very thick] (0,0)--++(5,0)--++(0,5)--++(-5,0)--cycle;
                \clip (0,0)--++(5,0)--++(0,5)--++(-5,0)--cycle;
                \foreach \x in {0,1,2,3,4,5}{
                \draw[-] (\x,0) --++(0,5);
                \draw[-] (0,\x) --++(5,0);
    }
    \two 00
    \two 01
    \two 04
    \two 12
    \two 13
    \two 21
    \two 22
    \two 23
    \two 32
    \two 42
    \two 43
    
    \end{scope}
    \begin{scope}[xshift=3cm]
        \draw [very thick] (0,0)--++(5,0)--++(0,1)--++(-1,0)--++(0,1)--++(-1,0)--++(0,1)--++(-1,0)--++(0,1)--++(-1,0)--++(0,1)--++(-5, 0)--++(0,-1)--++(1,0)--++(0,-1)--++(1,0)--++(0,-1)--++(1,0)--++(0,-1)--++(1,0)-- (0,0);
        \clip (0,0)--++(5,0)--++(0,1)--++(-1,0)--++(0,1)--++(-1,0)--++(0,1)--++(-1,0)--++(0,1)--++(-1,0)--++(0,1)--++(-5, 0)--++(0,-1)--++(1,0)--++(0,-1)--++(1,0)--++(0,-1)--++(1,0)--++(0,-1)--++(1,0)-- (0,0);
                \foreach \x in {-4,-3,-2,-1,0,1,2,3,4,5,6,7,8,9}{
                \draw[-] (\x,0) --++(0,5);
                \draw[-] (0,\x) --++(5, 0);
                \draw[-] (0,\x) --++(-5, 0);
                }
    \end{scope}
    \two 30
    \two 32
    \two 41
    \two 42
    \two 43
    \two 52
    \two 21
    \two 22
    \two 23
    \two 13
    \two -4 7
    
\end{tikzpicture}
    \caption{Left: A $({\protect\sqmo}, 2)$-prismatic polyomino. Right: A $({\protect\zmo}, 2)$-prismatic polyomino, which is the  row shift of the polyomino on the left.}
    \label{fig:Zprismatic}
\end{figure}

\begin{theorem}\label{thm:sq_and_z}
    For any positive integer $n$, the row shift gives a bijection from  $(\sqmo, n)$-prismatic polyominoes to $(\zmo, n)$-prismatic polyominoes.
\end{theorem}

\begin{proof}
    Let $P$ is a $(\sqmo, n)$-prismatic polyomino. Then  $f(P)$ is a $(\zmo,n)$-dense polyomino by Lemma \ref{thm:dense_bij}, and its coloring is $(\zmo,n)$-de Bruijn by Lemma \ref{thm:color_bij}; therefore $f(P)$ is $(\zmo,n)$-prismatic. So the row shift maps $(\sqmo, n)$-prismatic polyominoes to $(\zmo,n)$-prismatic polyominoes. 
    Similarly, the inverse transformation maps $(\zmo, n)$-prismatic polyominoes to $(\sqmo,n)$-prismatic polyominoes.
    Because the transformation is invertible, it is bijective.
\end{proof}

An example of this transformation is shown in Figure 
    \ref{fig:Zprismatic}.

\begin{remark}
Some of what we know about $(\sqmo, n)$-pris\-mat\-ic poly\-om\-inoes applies to $(\zmo, n)$-pris\-mat\-ic poly\-om\-inoes.
There is a unique shape for $(\zmo, n)$-pris\-mat\-ic poly\-om\-inoes. There are exactly 800 differently colored $(\zmo, 2)$-prismatic polyominoes, and at least $n!^{n} \cdot (n^2)!$ many $(\zmo, n)$-prismatic polyominoes. 

\end{remark}

\subsection{A Bijection for T and L Tetrominoes}
In this section, we will discuss the relationship between T and L tetrominoes under the row shift.
The row shift of a T tetromino is always an L tetromino, and our results from Section \ref{sec:sq_and_z} have analogs for these shapes, which we state without proof.

\begin{lemma}
    Let $p$ and $q$ be sets of cells such that $f(p) = q$. Then $p$ is an instance of the T tetromino if and only if $q$ is an instance of the L tetromino.
\end{lemma}

\begin{lemma}
    Let $P$ and $Q$ be sets of cells such that $f(P)=Q$. Then for any positive integer $N$, $P$ is a $(\tmo,N)$-dense polyomino if and only if $Q$ is a $(\lmo,N)$-dense polyomino.
\end{lemma}

\begin{lemma}
    For any positive integer $n$, a set of cells $S$ with an $n$-coloring includes exactly one instance of each $n$-coloring of the T tetromino if and only if $f(S)$ includes exactly one instance of each $n$-coloring of L tetromino. 
\end{lemma}

\begin{theorem}\label{thm:prismatic_lmo}
    For any positive integer $n$, the row shift gives a bijection from  $(\tmo, n)$-prismatic polyominoes to $(\lmo, n)$-prismatic polyominoes.
\end{theorem}

\begin{figure}
    \centering
    \begin{tikzpicture}[scale=.5]
    \begin{scope}[xshift= -9cm]
    \draw [very thick] (0,0)--++(9,0)--++(0,1)--++(-1,0)--++(0,1)--++(-1,0)--++(0,1)--++(-1,0)--++    (0,1)--++(-1,0)--++(0,1)--++(-1,0)--++(0,-1)--++(-1,0)--++(0,-1)--++(-1,0)--++(0,-1)--++(-1,0)--++(0,-1)--++(-1,0)--++(0,-1)--cycle;
    \clip (0,0)--++(9,0)--++(0,1)--++(-1,0)--++(0,1)--++(-1,0)--++(0,1)--++(-1,0)--++    (0,1)--++(-1,0)--++(0,1)--++(-1,0)--++(0,-1)--++(-1,0)--++(0,-1)--++(-1,0)--++(0,-1)--++(-1,0)--++(0,-1)--++(-1,0)--++(0,-1)--cycle;
                \foreach \x in {0,1,2,3,4,5,6,7,8,9}{
                \draw[-] (\x,0) --++(0,5);
                \draw[-] (0,\x) --++(9,0);
    }
    \two 44
    \two 22
    \two 32
    \two 20
    \two 11
    \two 31
    \two 50
    \two 60
    \two 70
    \two 51
    \two 61
    \two 71
    \end{scope}
    
    \begin{scope}[xshift=3cm]
    \draw [very thick] (0,0)--++(9,0) --++(0,1)--++(-2,0) --++(0,1)--++(-2,0) --++(0,1)--++(-2,0) --++(0,1)--++(-2,0) --++(0,1)--++(-1,0) -- cycle;
    
    \clip (0,0)--++(9,0) --++(0,1)--++(-2,0) --++(0,1)--++(-2,0) --++(0,1)--++(-2,0) --++(0,1)--++(-2,0) --++(0,1)--++(-1,0) -- cycle;
                \foreach \x in {0,1,2,3,4,5,6,7,8,9}{
                \draw[-] (\x,0) --++(0,5);
                \draw[-] (0,\x) --++(9,0);
    }
    \two 20
    \two 50
    \two 60
    \two 70
    
    \two 01
    \two 21
    \two 41
    \two 51
    \two 61
    
    \two 02
    \two 12
    
    \two 04
    
    \end{scope}
\end{tikzpicture}
    \caption{Left: A $({\protect\tmo}, 2)$-prismatic polyomino. Right: A $({\protect\lmo}, 2)$-prismatic polyomino, which is the row shift of the polyomino on the left.}
    \label{fig:T_prismatic}
\end{figure}

\begin{remark}
    Recall, we have not shown that $(\tmo,n)$-prismatic polyominoes exist for every positive integer $n$. Theorem \ref{thm:prismatic_lmo} hold vacuously in these cases. We have shown there are 168 distinct $(\tmo, 2)$-prismatic polyominoes, so there are 168 distinct $(\lmo, 2)$-prismatic polyominoes as well. These are related by the row shift, as demonstrated in  Figure \ref{fig:T_prismatic}.
\end{remark}

\section{The L Tromino} \label{sec:Lmo}
    In this section we discuss the shapes of $(\tro,8)$-dense polyominoes, in order to find $(\tro,2)$-prismatic polyominoes.
    We will refer to the cells in the L tromino as the \textit{central}, \textit{right}, and \textit{top} cell as shown in Figure \ref{fig:L_tromino_cells}.

    \begin{figure}[h]
\centering
\begin{tikzpicture}[scale=1.25]
    \draw[-] (1,1)--++(2,0)--++(0,-1)--++(-2,0)--++(0,2)--++(1,0)--++(0,-2);
    \node (C) at (1.5,.5){central};
    \node (R) at (2.5,.5){right};
    \node (T) at (1.5,1.5){top};
\end{tikzpicture}
\caption{The cells of the L tromino.}
\label{fig:L_tromino_cells}
\end{figure}

    For any $(\tro,8)$-dense polyomino $P$, we let $X_c$ denote the set of those cells of $P$ that are the central cell of some instance of the L tromino within $P$. 
    
We let $\Sr$ (resp. $\St$) denote the set of cells that belong to exactly one instance of the L tromino in $P$ and are the right cell (resp. the top cell) of that instance.

We let $\Srt$ denote the set of cells that are the right cell of one instance of the L tromino, and the top cell of another instance of the L tromino, but are not the central cell of any instance of the L tromino.

    \begin{lemma}\label{thm:bounds}
        For any $(\tro,8)$-dense polyomino $P$ with width $w$ and height $h$,
        \begin{equation}\label{eq:BCD}
            |P| = |X_c| + |\Sr|  + |\St|+ |\Srt|,
        \end{equation}
        and
        \begin{equation}\label{eq:dim_bound}
            |P| \geq |X_c| + \max(w,h).
        \end{equation}
    \end{lemma}

    \begin{proof}
        As every cell in $P$ belongs to at least one instance of $\tro$ as a {central} cell, {right} cell, or {top} cell, and the sets $X_c$, $\Sr$, $\St$, and $\Srt$ are disjoint, these sets partition the set of cells of $P$; Equation \ref{eq:BCD} follows immediately.
        
        The right-most cell of any row of $P$ cannot be the central cell for any instance of $\tmo$, and so must be an element of $\St \cup \Sr \cup \Srt$. Thus,
        $$|\St \cup \Sr \cup \Srt| \geq h.$$
        Similarly, the top cell of any column of $P$ is an element of $\St \cup \Sr \cup \Srt$, so 
        $$|\St \cup \Sr \cup \Srt| \geq w.$$
        Therefore,
        $$|\St \cup \Sr \cup \Srt| \geq  \max(w,h).$$
        Equation \ref{eq:dim_bound} follows:
        $$|P| = |X_c| + |\St|  + |\Sr|+ |\Srt| = |X_c|+|\St \cup \Sr \cup \Srt| \geq |X_c|+ \max(w,h).$$
        \end{proof}

    We apply this lemma to get a lower bound on the size of $P$.
        \begin{lemma}\label{thm:k>11}
            Every $(\tro,8)$-dense polyomino $P$ has $|P| \geq 12$, and has either a width or height of at least 4.
        \end{lemma}

        \begin{proof}
            Since each instance of $\tro$ in $P$ has a distinct central cell, $|X_c| \geq 8$. Since all cells in the right column of $P$ are elements of $\Sr$, $|\Sr| \geq 1$. Since all cells in the top row of $P$ are elements of $\St$, $|\St| \geq 1$.
            Then $|P| \geq 8+1+1 = 10$. Thus, $P$ cannot fit within a $3 \times 3$ square, and must have a width or height of at least~4. Then $|P| \geq 8+4 = 12$ by Equation \ref{eq:dim_bound}.
        \end{proof}

    We now tighten the lower bound on the size of $P$.
    \begin{lemma}\label{thm:k>12}
        Every $(\tro,8)$-dense polyomino $P$ has $|P| \geq 13$.
    \end{lemma}

    \begin{proof}
        Suppose for the sake of contradiction that $P$ is a $(\tro,8)$-dense polyomino  with  $|P| \leq 12$. Let $P$ have width $w$ and height $h$. By Lemma \ref{thm:k>11} $P$ has exactly 12 cells, so by Equation \ref{eq:dim_bound} $P$ has width and height at most 4. Therefore $P$ fits inside a $4 \times 4$ rectangular polyomino, which we call $R$.  On the diagram below, we label each cell of $R$ with the number of instances of $\tro$ that cell belongs to.
        \begin{center}
        \begin{tikzpicture}[scale=.6]
            \foreach \x in {0,1,2,3,4}{
                \draw[-] (\x,0) --++(0,4);
                \draw[-] (0, \x) --++(4,0);
            } 
        \clabel {0}{3}{1} \clabel {1}{3}{1} \clabel {2}{3}{1} \clabel {3}{3}{0}
        \clabel {0}{2}{2} \clabel {1}{2}{3} \clabel {2}{2}{3} \clabel {3}{2}{1}
        \clabel {0}{1}{2} \clabel {1}{1}{3} \clabel {2}{1}{3} \clabel {3}{1}{1}
        \clabel {0}{0}{1} \clabel {1}{0}{2} \clabel {2}{0}{2} \clabel {3}{0}{1}
        \end{tikzpicture}
        \end{center}
        We can construct $P$ by removing four cells from $R$.
        As there are only nine instances of $\tro$ in $R$, and eight instances of $\tro$ in $P$, we cannot remove any cell that belongs to multiple instance of $\tro$. Among those cells labeled 1, only the two cells adjacent to the cell labeled 0 belong to the same instance of $\tro$; we can therefore remove at most two cells labeled with a 1. The only other cell we can remove is the cell labeled 0. Since we cannot remove more than three cells from $R$ to construct a polyomino with at least eight instances of $\tro$, we have arrived at a contradiction.
    \end{proof}

    This argument suggests a shape for $P$, which we will use to show $|P| = 13$, and that the width and height of $P$ are at most 5.
    
    \begin{theorem}
        \label{thm: k=13}
        Every $(\tro,8)$-dense polyomino $P$ has $|P| = 13$, and has width and height at most 5.
    \end{theorem}
    \begin{proof}
        Following Lemma \ref{thm:k>12}, to show $|P|=13$ it suffices to exhibit a polyomino that has eight instances of $\tro$ and 13 cells, which we show here.
    \begin{center}
        \polycolored{(0,0) --++(4,0)--++(0,2)--++(-1,0)--++(0,1)--++(-1,0)--++(0,1)--++(-2,0) -- cycle}{4}{}
    \end{center}
    Suppose $P$ is any $(\tro,8)$-dense polyomino, and has width $w$ and height $h$.
    By Equation \ref{eq:dim_bound}, since $P$ has 13 cells, $13 \geq 8 + \max(w,h)$, so $5 \geq w$ and $5 \geq h$.
    \end{proof}

\begin{remark}
    Theorem \ref{thm: k=13} shows that every $(\tro,8)$-dense polyomino fits within a $5 \times 5$ box, which facilitates a brute-force search for their shapes. We find there are nine such shapes. We exhibit these in Figure \ref{fig:trominoessolution}, with $(\tro,2)$-de Bruijn colorings; these are therefore $(\tro,2)$-prismatic.
\end{remark}

    \begin{figure}
        \centering
        \noindent
        \begin{tabular}{cc|ccc}
        &\polycolored{(1,0) --++(3,0)--++(0,2)--++(-1,0)--++(0,1)--++(-1,0)--++(0,1)--++(-1,0)--++(0,1)--++(-1,0)--++(0,-4)--++(1,0) -- cycle}{4}{\two 02 \two 10 \two 11 \two 21 \two 22 \two 31}
        &\polycolored{(0,0) --++(3,0)--++(0,1)--++(1,0)--++(0,1)--++(-1,0)--++(0,1)--++(-1,0)--++(0,1)--++(-1,0)--++(0,1)--++(-1,0)--++(0,-5) -- cycle}{4}{\two 11 \two 21 \two 22 \two 12 \two 03}
        &\polycolored{(0,0) --++(4,0)--++(0,1)--++(-1,0)--++(0,2)--++(-1,0)--++(0,1)--++(-1,0)--++(0,1)--++(-1,0)-- cycle}{4}{\two 01 \two 11 \two 12 \two 21 \two 20}
        &\polycolored{(1,0)--++(4,0)--++(0,1)--++(-1,0)--++(0,1)--++(-1,0)--++(0,1)--++(-1,0)--++(0,1)--++(-1,0)--++(0,1)--++(-1,0)--++(0,-3)--++(1,0)--cycle}{4}{\two 20 \two 21 \two 31 \two 11 \two 12}\\[10pt]
        \polycolored{(0,0) --++(4,0)--++(0,2)--++(-1,0)--++(0,1)--++(-1,0)--++(0,1)--++(-2,0) -- cycle}{4}{\two 20 \two 21 \two 31 \two 02 \two 12 \two 03}
        &\polycolored{(1,0)--++(4,0)--++(0,1)--++(-1,0)--++(0,1)--++(-1,0)--++(0,1)--++(-1,0)--++(0,1)--++(-2,0)--++(0,-3)--++(1,0)--cycle}{4}{\two 03 \two 10 \two 11 \two 12 \two 20 \two 22}
        &\polycolored{(0,0)--++(5,0)--++(0,1)--++(-1,0)--++(0,1)--++(-1,0)--++(0,1)--++(-1,0)--++(0,1)--++(-1,0)--++(0,-1)--++(-1,0)--cycle}{4}{\two 11 \two 21 \two 22 \two 12 \two 30}
        &\polycolored{(0,0)--++(5,0)--++(0,1)--++(-1,0)--++(0,1)--++(-1,0)--++(0,1)--++(-2,0)--++(0,1)--++(-1,0)--cycle}{4}{\two 10 \two 11 \two 21 \two 12 \two 02}
        &\polycolored{(2,0)--++(3,0)--++(0,1)--++(-1,0)--++(0,1)--++(-1,0)--++(0,1)--++(-1,0)--++(0,1)--++(-1,0)--++(0,1)--++(-1,0)--++(0,-4)--++(2,0)--cycle}{4}{\two 02 \two 12 \two 13 \two 11 \two 21}\\
        \end{tabular}
        \caption{There are nine shapes of $({\protect\tro},2)$-prismatic polyomino. A branch-and-prune search shows that the three shapes to the left of the vertical bar each have 8 $({\protect\tro},2)$-de Bruijn colorings, and the six shapes right of the vertical bar each have 28 $({\protect\tro},2)$-de Bruijn colorings.}
        \label{fig:trominoessolution}
    \end{figure}

    \section{Conclusion}
    
    While the tools of linear programming and polyomino transformations have been useful to us, we expect that more strategies will be needed to answer open questions about de Bruijn polyominoes and dense polyominoes. 
    We consider it noteworthy that for T trominoes there were several shapes of prismatic polyomino, while for tetrominoes all evidence points toward unique shapes.

    \section*{Acknowledgments}

    Many thanks to the organizers of the Polymath Jr program, and especially to Steve Miller for organizing the project group at Williams College. This work was conducted with support from NSF Grant DMS2218374. 

\printbibliography

@article{BE,
author = {T. van Aardenne-Ehrenfest and N. G. de Brujin},
title={Circuits and Trees in Oriented Linear Graphs},
journal={Simon Stevin},
year = {1951},
volume={28},
pages={203-217}
}

@article{Cock,
    author = {Cock, J. C.},
    title = {Toroidal tilings from de Bruijn-Good cyclic sequences},
    journal = {Discrete Math},
    year = {1988},
    volume = {70},
    pages = {209-210}
}

@misc{Dow,
  author = {Mark Dow},
  title = {Minimal arrays containing all sub-array combinations of symbols: De Bruijn sequences and tori},
  year = {2012}, 
  howpublished = {\url{https://web.archive.org/web/20120427194606/http://lcni.uoregon.edu/~dow/Geek_art/Minimal_combinatorics/Minimal_arrays_containing_all_combinations.html}},
  note = {Online: accessed August 11 2023}
}

@book{Golomb,
    author = {Solomon W. Golomb},
    title = {Polyominoes: Puzzles, Patterns, Problems, and Packings - Revised and Expanded Second Edition},
    publisher = {Princeton University Press},
    year = {1994}
}

@article{DV,
author={Dale E. Varberg},
journal={The American Mathematical Monthly},
title={Pick's Theorem Revisited},
month={8},
year={1985},
volume={92},
number={8},
pages={584-587}
}

\appendix

\newpage

\section{Solutions to Problem 1.1} \label{app:16sols}

The following Python code generates all solutions to Problem \ref{prob:sixteen_squares}. We encode a 2-coloring of a $5 \times 5$ square polyomino $P$ as a binary string
$$s_{1,1}s_{1,2}s_{1,3}s_{1,4}s_{1,5}s_{2,1}s_{2,2}\dots s_{5,4}s_{5,5}$$
where $s_{i,j}$ represents the color of the cell in the $i$th row and $j$th column of $P$.

We perform a branch-and-prune search for strings representing partial colorings of $P$ in which no two instances of $\sqmo$ are colored the same way. 
We start with a list of the 32 binary strings of length 5, each string representing a coloring of the first row of $P$. 
Then we list all of the ways that five more bits may be appended to the end of each of these strings, coloring the next row of $P$. If a partial coloring includes multiple instances of $\sqmo$ colored the same way, we remove its string from the list. We iterate this process, coloring one row at a time. Throughout the process, a coloring of the first $r$ rows is represented as a string of $5r$ characters,
$$s_{1,1}s_{1,2}s_{1,3}s_{1,4}s_{1,5}s_{2,1}s_{2,2}\dots s_{r,4}s_{r,5}.$$

The following function evaluates a string representing a partial coloring, and returns \texttt{True} if no two instances of $\sqmo$ are colored the same way.

\begin{lstlisting}
def no_repeats(string): 
    squarelist = [] 
    for i in range(len(string)-6):
        if (i%5 < 4):
            squarelist.append(
                string[i] + string[i+1] + 
                string[i+5] + string[i+6])
    return len(set(squarelist)) == len(squarelist) 
\end{lstlisting}

The following function builds colorings one row at a time, and purges any that include two instances of $\sqmo$ that are colored the same way. The function returns a list of strings representing all 800 $(\sqmo,2)$-prismatic polyominoes, with the colors $\{0,1\}$ instead of $\{1,2\}$.

\begin{lstlisting}
def prismatic(): 
    rows = [''] 
    for i in range(5):
        rows = [string+symb for symb in '01' for string in rows]
    
    gridlist = rows
    for i in range(4): 
        gridlist = [grid+row for row in rows for grid in gridlist] 
        gridlist = [grid for grid in gridlist if no_repeats(grid)] 
    return gridlist
\end{lstlisting}

Variations on this code can be used to construct prismatic polyominoes with 2 colors for $\zmo$, $\lmo$, $\tmo$ and $\tro$. 
\end{document}